\documentclass[letterpaper,12pt,twoside,reqno]{amsart}
\usepackage[bottom=20mm,top=20mm]{geometry}
\usepackage{times}
\usepackage{amsmath}
\usepackage{amssymb}
\usepackage{amsfonts}
\usepackage{amsthm}

\theoremstyle{plain}
\newtheorem{theorem}{Theorem}

\newtheorem*{proposition*}{Proposition}

\newtheorem*{corollary*}{Corollary}

\newtheorem*{theorem*}{Theorem}
\newtheorem*{lemma*}{Lemma}
\newtheorem*{conjecture*}{Conjecture}

\newtheorem*{question*}{Question}

\theoremstyle{definition}

\newtheorem*{exercise*}{Exercise}

\theoremstyle{remark}

\newtheorem*{remark*}{Remark}
\newtheorem{remsTh}[theorem]{Remarks}

\newcommand{\subclass}[1]{}

\newcommand{\enumTi}[1]{\renewcommand{\theenumi}{#1}}

\newcommand{\alphenumi}{\enumTi{\alph{enumi}}}

\newcommand{\romenumi}{\enumTi{\roman{enumi}}}

\renewcommand{\em}{\sl}

\DeclareMathOperator{\conv}{conv}

\newcommand{\polar}{\vartriangle}

\newcommand{\RR}{\mathbb{R}}

\newcommand{\ZZ}{\mathbb{Z}}

\newcommand{\iprod}{\bullet}

\newlength{\algotabbingwidth}
\renewcommand{\iprod}{\cdot}

\begin{document}

\title[Traveling Salesman Polyhedra and Metric Cone]{%
  A note on the relationship between the Graphical Traveling Salesman Polyhedron, the Symmetric Traveling Salesman Polytope, and the Metric Cone\\
  {\footnotesize\textrm{(short communication)}}
}%
\author{Dirk Oliver Theis}%
\address{Dirk Oliver Theis, %
  Service de G\'eom\'etrie Combinatoire et Th\'eorie des Groupes, %
  D\'epartement de Math\'ematique, %
  Universit\'e Libre de Bruxelles, %
  Brussels, Belgium}%
\email{Dirk.Theis@ulb.ac.be}%
\thanks{Research supported by \textit{Deutsche Forschungsgemeinschaft (DFG)} as project RE~776/9-1.
  Author currently supported by \textit{Communaut\'e fran\c caise de Belgique -- Actions de Recherche
    Concert\'ees.}}

\date{Sat Feb 20 16:04:36 CET 2010}

\newcommand{\Vset}[1]{V_{#1}}


\begin{abstract}
  In this short communication, we observe that the Graphical Traveling Salesman Polyhedron is the intersection of the positive orthant with the Minkowski sum of
  the Symmetric Traveling Salesman Polytope and the polar of the metric cone.
  This follows almost trivially from known facts.  There are two reasons why we find this observation worth communicating none-the-less: It is very surprising;
  it helps to understand the relationship between these two important families of polyhedra.
\end{abstract}
\maketitle

\newcommand{\Esetn}{{E_n}}

\section{Introduction}

The \textit{Symmetric Traveling Salesman Polytope} is the convex hull of all characteristic vectors of edge sets of cycles (i.e., circuits) on the vertex set
$\Vset{n}:=\{1,\dots,n\}$ (in other words, Hamiltonian cycles in the complete graph with vertex set $\Vset{n}$).  For the formal definition, denote by $\Esetn$
the set of all two-element subsets of $\Vset{n}$.  This is the set of all possible edges of a graph with vertex set $\Vset{n}$.  The Symmetric Traveling
Salesman Polytope is then the following set:
\begin{equation*}
  S_n := \conv \Bigl\{ \chi^{C} \mid C \text{ is the edge set of a Hamiltonian cycle with vertex set $\Vset{n}$} \Bigr\} \subset \RR^\Esetn.
\end{equation*}
Here, for an edge set $F$, $\chi^F$ is the characteristic vector in $\RR^\Esetn$ with $\chi^F_e=1$ if $e\in F$, and zero otherwise.
The importance of the Symmetric Traveling Salesman Polytope comes mainly, but not exclusively, from its use in the solution of the so-called Symmetric Traveling
Salesman Problem, which consists in finding a Hamiltonian cycle of minimum cost.

The \textit{Graphical Traveling Salesman Polyhedron} is the convex hull of all characteristic vectors of edge multi-sets of connected Eulerian multi-graphs on
the vertex set $\Vset{n}$.  A multi-graph with vertex set $\Vset{n}$ has as its edge set a sub-multi-set of $\Esetn$, which is to say that our multi-graphs can
have parallel edges but no loops.
By defining, for any multi-set $F$ of edges of $K_n$, its characteristic vector $\chi^F \in \RR^\Esetn$ in such a way that $\chi^F_e$ counts the number of
occurrences of $e$ in $F$,
the Graphical Traveling Salesman Polyhedron is formally defined as 
\begin{multline*}
  P_n := \conv\Bigl\{ \chi^F \bigm| F \text{ is the edge multi-set of a connected Eulerian multi-graph}\\
  \text{with vertex set $\Vset{n}$} \Bigr\} \subset \RR^\Esetn.
\end{multline*}

Ever since the seminal work of Naddef \& Rinaldi \cite{NadRina91,NadRina93} on the two polyhedra, $P_n$ is considered to be an important tool for investigating
the facets of $S_n$.  Moreover, in works of Carr \cite{Carr2004} and Applegate, Bixby, Chv\`atal \& Cook \cite{ABCC2001}, $P_n$ has been used algorithmically in
contributing to solution schemes for the Symmetric Traveling Salesman Problem.

Numerous authors have expressed how close the connection between Graphical and Symmetric Traveling Salesman Polyhedra is.  The most basic justification for this
opinion is the fact that $S_n$ is a face of $P_n$ --- consisting of all points $x$ whose ``degree'' is two at every vertex: $\sum_{v\ne u} x_{uv} = 2$ for all
$u\in\Vset{n}$.
However, the connections are far deeper (see \cite{Gutin_Naddef} or \cite{OsRlTheis07} and the references therein).
In this short communication, we contribute the following surprising geometric observation to the issue of the relationship between these two polyhedra:

\begin{theorem*}
  $P_n$ is the intersection of the positive orthant with the Minkowski sum of $S_n$ and the polar $C_n^\polar$ of the metric cone $C_n$:
  \begin{equation}\label{eq:main-statement}
    P_n = (S_n+C_n^\polar)\cap \RR_+^\Esetn
  \end{equation}
\end{theorem*}

The metric cone consists of all $a\in\RR^\Esetn$ which satisfy the triangle inequality:
\begin{equation}\label{eq:triangle-ieq}
  a_{uv}\le a_{uw} + a_{wv}
\end{equation}
for all pairwise distinct vertices $u,v,w \in \Vset{n}$.  Consequently, its polar is generated as a cone by the vectors (we abbreviate $\chi^{\{e\}}$ to $\chi^e$)
\begin{equation}\label{eq:shortcut}
  \chi^{uw} + \chi^{wv} - \chi^{uv}.
\end{equation}

The proof of this theorem is an application of three or four known facts or techniques in the area of Symmetric and Graphical Traveling Salesman polyhedra.


\section{Proof}

We start with showing that $P_n\subset (S_n+C_n^\polar)\cap \RR_+^\Esetn$.  While $P_n\subset \RR_+^\Esetn$ holds trivially, $P_n\subset S_n+C_n^\polar$
follows from an argument of \cite{NadRina93}, which we reproduce here for the sake of completeness.

Let $x\in\ZZ_+^\Esetn$ be a the characteristic vector of the edge multi-set of a connected Eulerian multi-graph $G$ with vertex set $\Vset{n}$.  We prove by
induction on the number $m$ of edges of $G$, that $x$ can be written as a sum of a cycle and a number of vectors~\eqref{eq:shortcut}.  If $m=n$, then there is
nothing to prove.  Let $m\ge n+1$.  There exists a vertex $w$ of degree at least four in $G$.
We distinguish two cases.  The easy case occurs when $G\setminus w$ is still connected.  Here, we let $u$ and $v$ be two arbitrary (possibly identical)
neighbors of $w$.  By either replacing the edges $uw$ and $wv$ of $G$ with the new edge $uv$, if $u\ne v$, or deleting $uw$ and $wv$, if $u=v$, one obtains a
connected Eulerian multi-graph $G'$ with fewer edges than $G$.  The change in the vector $x$ amounts to subtracting the expression~\eqref{eq:shortcut}: $x' =
x-(\chi^{uw} + \chi^{wv} - \chi^{uv})$, if $u\ne v$, and $x' = x-(\chi^{uw} + \chi^{wv}$, if $u=v$.
In the slightly more difficult case when the graph $G\setminus w$ has at least two connected components, we can let $u$ and $v$ be two neighbors of $w$ in
distinct components of $G\setminus w$.  This makes sure that the graph $G'$ is still connected.
We conclude by induction that $x'$, and hence $x$, can be written as a sum of a cycle and a number of vectors~\eqref{eq:shortcut}.

We now prove $P_n \supset (S_n+C_n^\polar)\cap \RR_+^\Esetn$.  For this, we show that any inequality which is facet-defining for $P_n$ is valid for
$(S_n+C_n^\polar)\cap \RR_+^\Esetn$.

We again invoke an argument from \cite{NadRina93}: Naddef \& Rinaldi have shown\footnote{%
  In fact, Proposition~2.2 of \cite{NadRina93} states that the facet-defining inequalities for $P_n$ fall into three classes --- one of which is the class of
  non-negativity inequalities and the other two satisfy the triangle inequality.  } %
that the inequalities defining facets of $P_n$ fall into one of two categories: the non-negativity inequalities $x_e \ge 0$, with $e\in\Esetn$ (or positive
scalar multiples thereof), or inequalities whose coefficient vectors satisfy the triangle inequality~\eqref{eq:triangle-ieq}.  We reproduce the proof of this
statement.

First recall that an inequality $a\iprod x \ge \alpha$ is said to be \textit{dominated} by another inequality $b\iprod x \ge \beta$, if the face defined by the
first inequality is contained in the face defined by the second inequality.

Suppose that $a\iprod x \ge \alpha$ is not dominated by a non-negativity inequality (it need not be define a facet, though), and let $u,v,w$ be three distinct
vertices in $\Vset{n}$.  Then there exists an $x\in\ZZ_+^\Esetn$ defining the edge multi-set of a connected Eulerian multi-graph $G$ which has an edge between
$u$ and $v$, such that $a\iprod x = \alpha$.  If we replace the edge $uv$ of $G$ by the two edges $uw$ and $wv$, then we obtain a connected Eulerian
multi-graph, whose edge multi-set is given, in terms of its characteristic vector, by $x' := x+\chi^{uw} + \chi^{wv} - \chi^{uv}$.  Now $a\iprod x' \ge \alpha$,
implies $a_{uw} + a_{wv} - a_{uv} \ge 0$, i.e., the triangle inequality.

We now conclude the proof of the inclusion $P_n \supset (S_n+C_n^\polar)\cap \RR_+^\Esetn$.  Let $a\iprod x \ge \alpha$ be an inequality which is facet-defining
for $P_n$.
First note that the non-negativity inequalities are clearly satisfied by the right hand side of~\eqref{eq:main-statement}.
Hence, using what we have just discussed, let us assume that $a$ satisfies the triangle inequality.  This means that $a$ is a member of the metric cone $C_n$.
Consequently, the inequality $a\iprod x \ge 0$ is valid for $C_n^\polar$.  Further, since $S_n\subset P_n$, the inequality $a\iprod x\ge \alpha$ is clearly
valid for $S_n$.  Hence the inequality is valid for $S_n+C_n^\polar$.

This concludes the proof of the theorem.

\bigskip\noindent%
Note that, en passant, we have proved the following.  If we define $P'_n$ to be the set of all $y\in\RR^\Esetn$ which satisfy $a\iprod y\ge \alpha$ for every
inequality $a\iprod x \ge \alpha$ defining a facet of $P_n$ but not being a scalar multiple of a non-negativity inequality, then we have $S_n+C_n^\polar \subset
P'_n$.

\section*{Acknowledgments}

Thanks are extended to the \textit{Deutsche Forschungsgemeinschaft}, DFG, for funding this research, and to the \textit{Communaut\'e fran\c caise de Belgique --
  Actions de Recherche Concert\'ees} for supporting the author during the time the paper was written.

\end{document}